# *A product representation for π*

**S R Holcombe**

Email: srholcombe@gmail.com

## *Abstract*

This short note delivers, via elementary calculations, the product representation,

$$\pi = e^{3/2} \prod_{n=2}^{\infty} \left(1 - \frac{1}{n^2}\right)^{n^2} e \, .$$

## *1. A product containing ζ(3).*

The derivation of the final product can be achieved in a variety of ways, and in fact many of the results in this note are obtainable via more general results in [2]. The final result obtained here, however, is absent in that paper. The author is not sure if the resulting product is novel, however the form of the product appears to be absent from a more than cursory glance through the most common texts containing series, products and constants (this includes the work of Euler and notebooks of Ramanujan along with the compendiums of Gradshteyn, Prudnikov, Jolley, and Finch, see [1,3-7]).

Here the derivation will follow its original path using elementary means arising from examining miscellaneous problems involving Fourier expansions over a square plate or well (elasticity, conductivity, oscillating well). In such physical systems it is not uncommon to encounter integrals and Fourier series similar to, or of the form,

$$S = \sum_{n=1}^{\infty} \frac{1}{n} \int_0^{\pi} \int_0^{\pi} \sin(nx)\sin(ny)\,dxdy = \sum_{n=1}^{\infty} \frac{\left((-1)^n - 1\right)^2}{n^3} = \frac{7}{2}\zeta(3) , \qquad (1.1)$$

where $\zeta(z)$ is the Riemann zeta function [9]. It is then not difficult to show

$$S = \mathrm{Re}(\bar{S}), \quad \text{where} \quad \bar{S} = \frac{1}{2} \int_0^{\pi} \int_0^{\pi} \log\left(\frac{\sin\left(\frac{1}{2}(x+y)\right)}{\sin\left(\frac{1}{2}(x-y)\right)}\right) dxdy \, . \qquad (1.2)$$

Rigour here has been omitted for brevity, but all functions and manipulations may be justified using well known results in [8,9]. Using the product representation of the sine function $\sin(x) = x \prod_{n=1}^{\infty} \left(1 - \frac{x^2}{n^2\pi^2}\right)$, one may show after some elementary though tedious calculation,

$$S = 15\pi^2\log(2) - 9\pi^2\log(3)$$
$$+\pi^2\sum_{n=2}^{\infty}\left\{n^2\log\left(\frac{(n^2-1)2^8n^6}{(4n^2-1)^4}\right)+2n\log\left(\frac{(1+n)(2n-1)^2}{(n-1)(1+2n)^2}\right)+\log\left(1-\frac{3}{4n^2-1}\right)\right\}. \quad (1.3)$$

Two of these series may be determined explicitly, the last by using known product representations of trigonometric functions. The second last series is not so obvious and requires a bit more work. Observe firstly, for example, that

$$\prod_{n=1}^{N}(1+n)^n = \Gamma^N(N+2)\left\{\prod_{n=1}^{N}\Gamma(n+1)\right\}^{-1}.$$

Using similar representations for the other terms in the second last series it may be shown

$$\prod_{n=2}^{N}\frac{(1+n)^{2n}(2n-1)^{4n}}{(n-1)^{2n}(1+2n)^{4n}} = \frac{81}{4\pi^2}\frac{(N+2)^{2N}(N+1)^{2N}}{(N+\frac{3}{2})^{4N}}\frac{\Gamma^4(N+\frac{1}{2})}{\Gamma^2(N)\Gamma^2(N+1)}.$$

Taking the limit as $N\to\infty$ yields

$$\pi = \frac{9}{2}\prod_{n=2}^{\infty}\frac{(n-1)^n(2n+1)^{2n}}{(n+1)^n(2n-1)^{2n}}.$$

Therefore

$$S = 15\pi^2\log(2) - 9\pi^2\log(3) - \pi^2\log\left(\frac{64\pi}{243}\right) + \pi^2\sum_{n=2}^{\infty}n^2\log\left(\frac{\left(1-\frac{1}{n^2}\right)}{\left(1-\frac{1}{4n^2}\right)^4}\right), \quad (1.4)$$

and from (1.1) the following known result is obtained [4],

$$\frac{9^2}{2^9}\pi\exp\left(\frac{7}{2\pi^2}\zeta(3)\right) = \prod_{n=2}^{\infty}\frac{\left(1-\frac{1}{n^2}\right)^{n^2}}{\left(1-\frac{1}{4n^2}\right)^{4n^2}}. \quad (1.5)$$

## 2. *A product for π*

Consider $a, x, \alpha_N \in \mathbb{R} : x > 1/4$ and the product,

$$\prod_{n=2}^{\infty}\left(1-\frac{1}{xn^2}\right)^{xn^2} a = \prod_{n=2}^{\infty}\exp\left(-1+O\left(1/xn^2\right)+\log(a)\right) \simeq \alpha_N\prod_{N>>2}^{\infty}\exp\left(-1+\log(a)\right).$$

For $a<e$ the product converges to zero, while for $a>e$ the product diverges. Hence consider now the following finite products

$$P_N = Q_N / U_N, \quad \text{where} \quad Q_N = \prod_{n=2}^{N}\left(1-\frac{1}{n^2}\right)^{n^2}, \quad U_N = \prod_{n=2}^{N}\left(1-\frac{1}{4n^2}\right)^{4n^2}. \tag{2.1}$$

Each product $Q_N$, $U_N$ converges to zero for increasing N, however their ratio does not necessarily converge to zero. The form of the product obtained in (1.5) therefore warrants a more careful inspection which in turn yields an expedient derivation of (1.5), or alternatively, a derivation of Euler's result [3]. Consider then,

$$P(x,y) = \prod_{n=2}^{\infty} \frac{\left(1-\dfrac{1}{xn^2}\right)^{xn^2}}{\left(1-\dfrac{1}{yn^2}\right)^{yn^2}} = A(x,y)\prod_{n=2}^{\infty}\left(1-\frac{1}{h(x,y)n^2}\right)^{h(x,y)n^2} e \tag{2.2}$$

where both $x,y>\frac{1}{4}$ and A and h are some as yet unspecified real functions of x and y. It may be shown that

$$\sum_{k=1}^{\infty}\frac{\left(\zeta(2k)-1\right)}{(k+1)}\left(\frac{1}{y^k}-\frac{1}{x^k}+\frac{1}{h(x,y)^k}\right) = \log A(x,y). \tag{2.3}$$

If it is specified that $h(x,y)=x$, then

$$\sum_{k=1}^{\infty}\frac{\left(\zeta(2k)-1\right)}{(k+1)y^k} = \log A(y), \qquad \text{or} \qquad A(y) = \frac{1}{\prod_{n=2}^{\infty}\left(1-\dfrac{1}{yn^2}\right)^{yn^2} e}. \tag{2.4}$$

Therefore

$$\prod_{n=2}^{\infty}\frac{\left(1-\dfrac{1}{xn^2}\right)^{xn^2}}{\left(1-\dfrac{1}{yn^2}\right)^{yn^2}} = \frac{\prod_{n=2}^{\infty}\left(1-\dfrac{1}{xn^2}\right)^{xn^2} e}{\prod_{n=2}^{\infty}\left(1-\dfrac{1}{yn^2}\right)^{yn^2} e}. \tag{2.5}$$

Now from the previously utilised product representation of sine one may readily establish

$$\log\left(\sin(\pi x)\right) = \log(\pi x) - \sum_{k=1}^{\infty}\frac{x^{2k}\zeta(2k)}{k}, \tag{2.6}$$

from which it follows

$$\sum_{k=1}^{\infty}\frac{x^{2k}\zeta(2k)}{k+1} = \frac{1}{2} - \log\left(\sin(\pi x)\right) + \frac{2}{x^2}\int_0^x t\log\left(\sin(\pi t)\right)dt. \tag{2.7}$$

Hence

$$A(y)=\frac{1}{\sin\left(\frac{\pi}{\sqrt{y}}\right)}\left(1-\frac{1}{y}\right)^{y}\exp\left\{\frac{3}{2}+R(y)\right\}, \tag{2.8}$$

where,

$$R(y)=2y\int_{0}^{1/\sqrt{y}} t\log\left(\sin(\pi t)\right)dt\,. \tag{2.9}$$

The general product of (2.2) may then be given as

$$P(x,y)=\frac{A(y)}{A(x)}. \tag{2.10}$$

Note that $R(1)=-\log(2)$ which is obtained via an elementary use of trigonometric identities. Noting also the easily determined limit

$$\lim_{x\to 1}\left(1-\frac{1}{x}\right)^{-x}\sin\left(\frac{\pi}{\sqrt{x}}\right)=\frac{\pi}{2}, \tag{2.11}$$

expressions (2.10) and (1.5) provide another derivation of Euler's result [3], i.e. $R(4)=-\log(2)+\frac{7}{2\pi^{2}}\zeta(3)$, while the explicit form of $1/A(1)$ yields the elegant but slowly converging product

$$\pi=e^{3/2}\prod_{n=2}^{\infty}\left(1-\frac{1}{n^{2}}\right)^{n^{2}}e\,. \tag{2.12}$$

## 3. Acknowledgements

The author would like to thank Steven Finch for helpful comments and references.

## 4. Addendum

Steven Finch recently (2020) published "Errata and Addenda to Mathematical Constants" on the arxiv.org [10]. In this addendum he has included a product which he originally suggested to me and to which I conveyed a derivation via private correspondence though did not originally include in this article. Since his errata and addenda essentially references this article, for completeness I have now included the derivation of this second and related infinite product. Namely consider

$$P(z)=\frac{\sin(\pi z)}{\left(1-z^{2}\right)^{\frac{1}{z^{2}}}}\exp\left\{-\frac{3}{2}-2\int_{0}^{1}t\log\left(\sin(\pi z t)\right)dt\right\} \tag{4.1}$$

Using the infinite product representation of sine one may readily show

$$\left(1-z^2\right)^{\frac{1}{z^2}} eP\left(z\right)=\prod_{n=1}^{\infty}\left(1-\frac{z^2}{n^2}\right)^{\frac{n^2}{z^2}} e\,, \tag{4.2}$$

where upon letting $z=i$ a similar product to (2.12) can be evaluated as

$$\prod_{n=1}^{\infty}\left(1+\frac{1}{n^2}\right)^{n^2}\frac{1}{e}=\frac{2}{eP\left(i\right)}\ . \tag{4.3}$$

Now note

$$S=\int_0^1 t\log\left(i\sinh\left(\pi t\right)\right)dt=\frac{i\pi}{4}-\frac{1}{2}\log\left(2\right)+\frac{\pi}{3}+\int_0^1 t\log\left(1-e^{-2\pi t}\right)dt \tag{4.4}$$

and

$$\begin{aligned}\int_0^1 t\log\left(1-e^{-2\pi t}\right)dt&=-\int_0^1 t\sum_{n=1}^{\infty}\frac{e^{-2\pi nt}}{n}dt=-\sum_{n=1}^{\infty}\frac{1}{n}\times\frac{1-\left(1+2\pi n\right)e^{-2\pi n}}{4\pi^2n^2}\\&=\sum_{n=1}^{\infty}\frac{\left(1+2\pi n\right)e^{-2\pi n}}{4\pi^2n^3}-\frac{1}{4\pi^2}\zeta\left(3\right)\end{aligned} \tag{4.5}$$

Hence

$$P\left(i\right)=4\sinh\left(\pi\right)\exp\left\{-\frac{3}{2}-\frac{2\pi}{3}+\frac{1}{2\pi^2}\zeta\left(3\right)-\sum_{n=1}^{\infty}\frac{\left(1+2\pi n\right)e^{-2\pi n}}{2\pi^2n^3}\right\} \tag{4.6}$$

By the definition of the polylogarithm Li, and use of (4.3) one has

$$\prod_{n=1}^{\infty}\left(1+\frac{1}{n^2}\right)^{n^2}\frac{1}{e}=\frac{1}{2\sinh\left(\pi\right)}\exp\left\{\frac{1}{2}+\frac{2\pi}{3}-\frac{1}{2\pi^2}\zeta\left(3\right)+\frac{1}{2\pi^2}\mathrm{Li}_3\left(e^{-2\pi}\right)+\frac{1}{\pi}\mathrm{Li}_2\left(e^{-2\pi}\right)\right\} \tag{4.7}$$